\documentclass[11pt]{article}
\usepackage{latexsym}
\usepackage[all]{xy}
\usepackage{amsmath}
\usepackage{amssymb}
\usepackage{amsfonts}
\usepackage[applemac]{inputenc}

\begin{document}

\title{\textbf{A note on the cotangent complex in derived algebraic geometry}}
\bigskip
\bigskip

\author{Gabriele Vezzosi\\
\small{Dipartimento di Matematica Applicata} \\ 
\small{Universit\`a di Firenze,}\\
\small{Firenze, Italy}}

\date{July 2010}

\maketitle

This short and elementary note originated from questions asked by many algebraic geometers when they were first exposed to a sketch of derived algebraic geometry. It is aimed at giving an example of how deformation theory takes a very transparent and uniform geometrical form when viewed into the world of derived algebraic geometry. In particular, we show how the cotangent complex of a scheme $X$ already contains informations about $X$ viewed as a derived scheme so that, in some sense, the cotangent complex was already, ``secretly'', an object in derived algebraic geometry long before this theory existed. By reversing this point of view, one could also say that, believing that the \emph{entire} cotangent complex is a natural geometric object leads one directly to the basic set-up of derived algebraic geometry. \\The proofs of the statements below are all in one form or another already contained in \cite{hagII}; the main point here is merely to extract or give a guide to the few needed results from \cite{hagII}, and put them in a form that could hopefully constitute a more direct answer to the various questions I have been asked.\\ I would like to thank Daniel Huybrechts for his questions and comments, Luchezar L. Avramov for useful explanations, and Bertrand To\"en for many recent conversations on this and related topics.\\

The few concepts and definitions we need from the theory of derived algebraic geometry, can be found in \cite{hagII}, or in the reviews \cite{hagdag, seattle}; another very good reference, though in a different framework, is \cite{lu}. As a tiny bit of notation, $\mathrm{Hom}_*$ will denote the Hom-sets in the category $\mathbf{Sch}_{*}$ of pointed schemes, while $\mathbb{R}\mathrm{Hom}_*$ will denote the Hom-sets in the homotopy category $\mathbf{dSch}_{*}$ of pointed derived schemes (or stacks). $\mathbf{dSt}_{k}$ will be the model category of derived stacks for the étale topology (denoted as $k-D^{-}Aff^{\sim ,\, \textrm{\'et}}$ in \cite[Def. 2.2.2.14]{hagII}). The acronym cdga stands for commutative differential graded algebra, with cohomological differential and concentrated in non-positive degrees.\\

\noindent \textbf{1. Derived interpretation of the cotangent complex --} Let $X$ be  a scheme over a field $k$ and $\mathbb{L}_{X}$ be its cotangent complex. Deformation theory of $X$ relates to the first two Ext-groups of $\mathbb{L}_{X}$ with values in $k$. If we pick a point $x$ in $X(k)$, the relevant spaces for classical deformation theory of $X$ around $x$ are $\mathrm{Ext}^{0}(\mathbb{L}_{X,x},k)$ and $\mathrm{Ext}^{1}(\mathbb{L}_{X,x},k)$, where we have put $\mathbb{L}_{X,x}:=x^{*}\mathbb{L}_{X}$. But while $\mathrm{Ext}^{0}(\mathbb{L}_{X,x},k)$ has a geometric and modular interpretation, as $\mathrm{Hom}_{*}(D_0:=\mathrm{Spec}\, k[\varepsilon_0], (X,x))$ (where $k[\varepsilon_0]$ is the usual ring of dual numbers over $k$), already $\mathrm{Ext}^{1}(\mathbb{L}_{X,x},k)$ does not have an analogous interpretation. More precisely, there is no pointed \emph{scheme} $D_1$ (abelian co-group object in $\mathbf{Sch}_{*,k}$) together with a natural (group) isomorphism  $$\mathrm{Ext}^{1}(\mathbb{L}_{X,x},k)\simeq \mathrm{Hom}_{*}(D_1, (X,x))$$ for each scheme $X$ and any $x\in X(k)$. Derived algebraic geometry solves this problem, and at the same time uniformly answers the same question for all $\mathrm{Ext}^{i}(\mathbb{L}_{X,x},k)$, $i\geq 0$: such a geometric/modular interpretation exists but one has to allow ``schemes'' with values in arbitrary cdga's or simplicial commutative algebras. $\mathbb{L}_{X}$ does encode (via its Ext-groups) precisely all \emph{deformations} of $X$, as long as one interprets these deformations as deformations in \emph{higher derived directions}. Let me make this more precise. And, in the hope of being more understandable, let me suppose that $k$ is of characteristic $0$, so that I can use cdga's (instead of simplicial commutative algebras) which are usually more familiar to algebraic geometers.\\

\noindent \textbf{Definition 1.1}  Let $i\in \mathbb{N}$ and $k[i]$ be the $k$-dg-module having $k$ just in degree $-i$; we can then consider the trivial square zero extension $k$-cdga $$ k[\varepsilon_i]:= k\oplus k[i].$$ Note that $k[\varepsilon_i]$ is concentrated in non-positive degrees (with a cohomological, i.e. degree-increasing, differential), and that its only nontrivial cohomology groups are concentrated in degrees $0$ an $-i$, where they both equal $k$. $k[\varepsilon_i]$ is called the \emph{dg-ring of $i$-th order dual numbers} over $k$, and its derived spectrum $\mathbf{D}_i:=\mathbb{R}\mathrm{Spec}\, k[\varepsilon_i]$ is called the \emph{derived  $i$-th order infinitesimal disk} over $k$. Note that, for odd $i$, $k[\varepsilon_i]$ is the free cdga on the $k$-dg module $k[i]$, i.e. on one generator in degree $-i$. \\

\noindent \textbf{Proposition 1.2.}  If $x\in X(k)$ is a $k$-rational point in $X$, then for each $i\in \mathbb{N}$, there is a canonical group isomorphism
$$\mathrm{Ext}_{k}^{i}(\mathbb{L}_{X,x},k)\simeq \mathbb{R}\mathrm{Hom}_{*}(\mathbf{D}_i, (X,x)),$$ where $\mathbb{R}\mathrm{Hom}_{*}=\mathrm{Hom}_{\mathrm{Ho}(\mathrm{Spec}\, k/dSt_{k})}.$\\

In other words, for any $i\in \mathbb{N}$, the functor
$$\mathbf{Sch}_{*,k}\longrightarrow \mathbf{Ab}\, \,\,: \, (X,x)\longmapsto \mathrm{Ext}_{k}^{i}(\mathbb{L}_{X,x},k),$$
while not co-representable in $\mathbf{Sch}_{*,k}$, it is indeed co-represented, by $\mathbf{D}_i$, in the larger (homotopy) category $\mathbf{dSch}_{*,k}$ of pointed derived schemes.
Therefore, this can also be used to give an alternative, more geometrical, definition of the full cotangent complex.\\

\noindent Proposition 1.2 is an immediate corollary of the following, more general\\

\noindent \textbf{Proposition 1.3.}  Let $R$ be a commutative $k$-algebra and $x:\mathrm{Spec}\, R\rightarrow X$ be a (classical) $R$-valued point of $X$. Then for each $i\in \mathbb{N}$, there is a canonical group isomorphism
$$\mathrm{Ext}_{R}^{i}(\mathbb{L}_{X,x},R)\simeq \mathbb{R}\mathrm{Hom}_{*}(\mathbf{D}_{i,R}, (X,x)),$$ where $\mathbb{R}\mathrm{Hom}_{*}:=\mathrm{Hom}_{\mathrm{Ho}(\mathrm{Spec}\, R/dSt_{k})},$ and $\mathbb{L}_{X,x}:=x^{*}\mathbb{L}_{X}$. \\

\noindent \textbf{Proof.} The proof is almost tautological, the only thing one really needs to remark is that Illusie's cotangent complex for a scheme $X$ coincide with the cotangent complex of $X$, viewed as a derived scheme or stack, as defined in \cite[\S 1.4.1.]{hagII}. Given this, we simply observe that $\mathrm{Ext}_{R}^{i}(\mathbb{L}_{X,x},R)\simeq \mathrm{Hom}_{\mathrm{D}(R)}(\mathbb{L}_{X,x},R[i])$ and, by \cite[1.4.1.5 (1), p. 99 and 1.4.1.4, p. 98 ]{hagII}  (with the notations used there)
$$\mathrm{Hom}_{\mathrm{D}(R)}(\mathbb{L}_{X,x},R[i])\simeq \pi_{0}(\mathbb{D}\mathrm{er}_X(R,R[i]))\simeq \mathrm{Hom}_{\mathrm{Ho}(\mathrm{Spec}\, R/dSt_{k})}(\mathbf{D}_{i,R}, (X,x)),$$ where $\mathbf{D}_{i,R}=\mathbb{R}\mathrm{Spec}(R)[R[i]]$, in the notations of \cite[p. 98]{hagII}, which is a co-group object.

 \hfill $\Box$ \\

\noindent \textbf{Remark 1.4.} One could also define, for any scheme $X$ over $k$ and any $i\in \mathbb{N}$, the $i$-th order derived tangent stack of $X$ as $$\mathbf{T}^{i}X:=\mathbb{R}\mathrm{HOM}_{\mathbf{dSt}_{k}}(\mathbf{D}_{i}, X),$$ where $\mathbb{R}\mathrm{HOM}_{\mathbf{dSt}_{k}}$ denotes the derived \emph{internal} Hom in $\mathbf{dSt}_{k}$ (so that $\mathbf{T}^{i}X$ is indeed a derived stack over $k$).  However the usefulness of this definition is doubtful, at least as far as the interpretation of the cotangent complex in derived algebraic geometry is concerned. In fact, already  $\mathbf{T}^{0}X\equiv \mathbf{T}X$ (in the notations of \cite[Def. 1.4.1.2]{hagII}) contains all the (cohomological) informations about $\mathbb{L}_{X}$.\\
%For example, there is a canonical isomorphism (in the derived category of $k$-modules) $$\mathbb{R}\mathrm{Hom}(\mathbb{L}_{X,x},k)\simeq \mathrm{DK}(\mathrm{Map}_{*}(\mathbf{D}_{0}, (X,x))),$$ where $\mathrm{Map}_{*}$ is the mapping space in the model category $\mathrm{Spec}\, k/dSt_{k}$ and $\mathrm{DK}: $ \\

\noindent \textbf{2. Global --} Of course, there is a corresponding \emph{global} version saying that, for any $i\in \mathbb{N}$, there is a canonical group isomorphism
$$\mathrm{Ext}_{\mathcal{O}_{X}}^{i}(\mathbb{L}_{X},\mathcal{O}_{X})\simeq \mathbb{R}\mathrm{Hom}_{*}(\mathbf{D}_{i,X}, X),$$
where $\mathbf{D}_{i,X}$ is the derived pullback of $\mathbf{D}_{i}$ to $X$, i.e. $\mathbf{D}_{i,X}=\mathbb{R}\underline{\mathrm{Spec}}_{X}(\mathcal{O}_{X}\oplus \mathcal{O}_{X}[i])$, and $\mathbb{R}\mathrm{Hom}_{*}$ now denotes the Hom-set in the homotopy category of $X/\mathbf{dSt}_{k}$. In other words, $\mathrm{Ext}_{\mathcal{O}_{X}}^{i}(\mathbb{L}_{X},\mathcal{O}_{X})$ is in bijection with \emph{retractions} of the canonical map $X\rightarrow \mathbf{D}_{i,X}$ (given by the first projection $\mathcal{O}_{X}\oplus \mathcal{O}_{X}[i]\rightarrow \mathcal{O}_{X}$). This follows from a straightforward analog of Proposition 1.3 :\\

\noindent \textbf{Proposition 2.1.}  \textbf{(i)} Let $R$  be a commutative $k$-algebra, $M$ a $R$-dg-module and $x:\mathrm{Spec}\, R\rightarrow X$ be a (classical) $R$-valued point of $X$. Then for each $i\in \mathbb{N}$, there is a canonical group isomorphism
$$\mathrm{Ext}_{R}^{i}(\mathbb{L}_{X,x},M)\simeq \mathbb{R}\mathrm{Hom}_{*}(\mathbf{D}_{i,R}(M), (X,x)),$$ where $\mathbb{R}\mathrm{Hom}_{*}:=\mathrm{Hom}_{\mathrm{Ho}(\mathrm{Spec}\, R/dSt_{k})},$  $\mathbb{L}_{X,x}:=x^{*}\mathbb{L}_{X}$ and  $\mathbf{D}_{i,R}(M):= \mathbb{R}\mathrm{Spec}(R\oplus M[i])$.\\
\indent \textbf{(ii)} Let $\mathcal{M}$ be a quasi-coherent $\mathcal{O}_{X}$-dg-Module on $X$. Then for each $i\in \mathbb{N}$, there is a canonical group isomorphism
$$\mathrm{Ext}_{\mathcal{O}_X}^{i}(\mathbb{L}_{X}, \mathcal{M})\simeq \mathbb{R}\mathrm{Hom}_{*}(\mathbf{D}_{i,X}(\mathcal{M}), X),$$ where $\mathbb{R}\mathrm{Hom}_{*}:=\mathrm{Hom}_{\mathrm{Ho}(X/dSt_{k})}$ and $\mathbf{D}_{i,X}(\mathcal{M}):=\mathbb{R}\underline{\mathrm{Spec}}_{X}(\mathcal{O}_{X}\oplus \mathcal{M}[i]).$ Note that $\mathbf{D}_{i,X}(\mathcal{M})$ is contravariantly functorial in $\mathcal{M}$. \\

\noindent \textbf{3. Kodaira-Spencer class and obstructions --} Using the geometric and modular interpretation of the cotangent complex given above, the Kodaira-Spencer class and obstructions are visible as \emph{morphisms} of derived schemes. Let's see this, to fix ideas, in the classical setup of a flat extension of a flat morphism through a square-zero closed immersion.
Let $$\xymatrix{X\ar[r]^-{f} & Y \ar[r] & S}$$ be morphisms of schemes, with $f$ flat, and $j:Y\hookrightarrow Y'$ be a closed immersion of $S$-schemes defined by a square-zero quasi-coherent Ideal $J$. It is well known (\cite[Th. 2.17 (ii)]{ill}) that there is an obstruction $\mathrm{Obs}(f;J)\in \mathrm{Ext}_{\mathcal{O}_{X}}^{2}(\mathbb{L}_{X/Y}, f^{*}J)$ which vanishes if and only if $f$ admits a lift to $f':X'\rightarrow Y'$ where $X\hookrightarrow X'$ is the closed immersion defined by the Ideal $f^{*}J$. Moreover, one has $$\mathrm{Obs}(f;J)= f^{*}[j] \times \mathrm{KS}(X/Y/S)$$ where $[j]\in \mathrm{Ext}_{\mathcal{O}_{Y}}^{1}(\mathbb{L}_{Y/S}, J)$ classifies $j$, $\mathrm{KS}(X/Y/S)\in \mathrm{Ext}_{\mathcal{O}_{X}}^{1}(\mathbb{L}_{X/Y}, f^{*}\mathbb{L}_{Y/S})$ is the Kodaira-Spencer class associated to $X\rightarrow Y \rightarrow S$ and $\times$ denotes the Yoneda composition product.  We will now recover both the Kodaira-Spencer and obstruction classes as morphisms between derived schemes.\\

\noindent - \textsf{Kodaira-Spencer morphism} - First let us notice that, by Proposition 2.1, the identity map $\mathbb{L}_{X/Y}\rightarrow \mathbb{L}_{X/Y}$ induces a morphism of derived stacks $$\mathbb{R}\underline{\mathrm{Spec}}_{X}(\mathcal{O}_{X}\oplus \mathbb{L}_{X/Y})\longrightarrow X$$ which is a retraction of the canonical map $X\rightarrow \mathbb{R}\underline{\mathrm{Spec}}_{X}(\mathcal{O}_{X}\oplus \mathbb{L}_{X/Y})$. Now, the transitivity triangle associated to $X\rightarrow Y \rightarrow S$ yields a map $\mathbb{L}_{X/Y}\rightarrow f^{*}\mathbb{L}_{Y/S}[1]$, which induces, by functoriality,   a map of derived schemes $$\mathbf{D}_{1,X}(f^{*}\mathbb{L}_{Y/S})=\mathbb{R}\underline{\mathrm{Spec}}_{X}(\mathcal{O}_{X}\oplus f^{*}\mathbb{L}_{Y/S}[1])\longrightarrow \mathbb{R}\underline{\mathrm{Spec}}_{X}(\mathcal{O}_{X}\oplus \mathbb{L}_{X/Y})=\mathbf{D}_{0,X}(\mathbb{L}_{X/Y}).$$ The composition (\emph{Kodaira-Spencer morphism}) $$\mathrm{ks}(X/Y/S): \, \mathbf{D}_{1,X}(f^{*}\mathbb{L}_{Y/S})\longrightarrow \mathbf{D}_{0,X}(\mathbb{L}_{X/Y})\longrightarrow X$$ gives back, by Proposition 2.1, an element $\mathrm{KS(X/Y/S)}\in \mathrm{Ext}_{\mathcal{O}_{X}}^{1}(\mathbb{L}_{X/Y}, f^{*}\mathbb{L}_{Y/S} )$: this is the classical Kodaira-Spencer class of $X\rightarrow Y \rightarrow S$. In a similar way one can get geometrically higher Kodaira-Spencer classes in higher cohomology.\\

\noindent - \textsf{Obstruction morphism} - Now we can recover the obstruction $\mathrm{Obs}(f;J)$ as a composite of morphisms between derived schemes. First observe that the shifted $j$-class $[j][1]\in \mathrm{Ext}_{\mathcal{O}_{Y}}^{2}(\mathbb{L}_{Y/S}, J)\simeq \mathrm{Ext}_{\mathcal{O}_{Y}}^{0}(\mathbb{L}_{Y/S}[1], J[2])$ is equivalent by Proposition 2.1, to a morphism $$\mathbf{D}_{2,X}(f^{*}J)=\mathbb{R}\underline{\mathrm{Spec}}_{X}(\mathcal{O}_{X}\oplus f^{*}J[2])\longrightarrow \mathbb{R}\underline{\mathrm{Spec}}_{X}(\mathcal{O}_{X}\oplus f^{*}\mathbb{L}_{Y/S}[1])=\mathbf{D}_{1,X}(f^{*}\mathbb{L}_{Y/S})$$ which can be composed with the Kodaira -Spencer morphism $\mathrm{ks}(X/Y/S)$ above to get a retraction (\emph{obstruction morphism}) $$\xymatrix{\mathrm{obs}(f;J): \mathbf{D}_{2,X}(f^{*}J)\ar[r] & \mathbf{D}_{1,X}(f^{*}\mathbb{L}_{Y/S})\ar[rr]^-{\mathrm{ks}(X/Y/S)} & & X}.$$ By Proposition 2.1, this morphism corresponds to a element in $\mathrm{Ext}_{\mathcal{O}_{X}}^{2}(\mathbb{L}_{X/Y}, f^{*}J)$ which is exactly the obstruction class $\mathrm{Obs}(f;J)$. Note that the Yoneda product obviously translates into a composition of morphisms of derived schemes.
\\

\noindent \textbf{4. Derived cotangent complex --} Moreover, the derived cotangent complex $\mathbb{L}_{F}$ (i.e. the cotangent complex as considered in derived algebraic geometry) is also defined, along the same lines, for any algebraic \emph{derived} (possibly higher) \emph{stack} $F$, and still have a similar geometric/modular interpretation which makes it natural and therefore computable.  Given any moduli underived $1$-stack $\mathcal{X}$, there is usually a nice \emph{derived moduli extension} $\mathcal{X}^{\mathrm{der}}$, i.e. a derived stack classifying an appropriate derived version of the geometric objects classified by $\mathcal{X}$, and whose truncation is $\mathcal{X}$. If one defines $\mathbb{L}^{\mathrm{der}}_{\mathcal{X}}:= \mathbb{L}_{\mathcal{X}^{\mathrm{der}}}$, then $\mathbb{L}^{\mathrm{der}}_{\mathcal{X}}$ is \emph{different}, in general, from the usual underived stacky cotangent complex of $\mathcal{X}$, which is often not known even for commonly studied moduli stacks; and it is $\mathbb{L}^{\mathrm{der}}_{\mathcal{X}}$ that controls the full deformation theory of the stack $F$. Even when the stacky cotangent complex agrees with a truncation of  $\mathbb{L}^{\mathrm{der}}_{\mathcal{X}}$, it is the latter that have nicer functorial properties and, as shown above in the case of schemes, a more natural geometrical interpretation. To give an example, let us consider $\mathcal{X}:=\mathbf{Vect}_{n}(X)$, the stack of rank $n$ vector bundles on a smooth projective variety $X$. There is a quite natural derived extension $\mathbf{Vect}_{n}(X)^{\mathrm{der}}$ classifying rank $n$ \emph{derived vector bundles} on $X$ ($\mathbf{Vect}_{n}(X)^{\mathrm{der}}$ is denoted as $\mathbb{R}\mathbf{Vect}_{n}(X)$ in \cite{hagII, hagdag}). And, if $x:\mathrm{Spec}\, k \rightarrow \mathbf{Vect}_{n}(X)^{\mathrm{der}}$ is a classical point, corresponding to a vector bundle $E\rightarrow X$, we have $$\mathbb{L}^{\mathrm{der}}_{\mathbf{Vect}_{n}(X),\, x}= \mathbb{L}_{\mathbf{Vect}_{n}(X)^{\mathrm{der}},\, x} \simeq \mathrm{C}_{\mathrm{Zar}}^{*}(X, \mathrm{End}(E))^{\vee}[-1]$$ which is exactly the complex everybody uses to compute the correct deformation theory of $\mathbf{Vect}_{n}(X)$ at $x$. But it \emph{cannot} be the stacky cotangent complex of \emph{any} algebraic $1$-stack (or even $n$-stack) locally of finite presentation, because it is perfect and of amplitude not in $[-1,\infty)$ for arbitrary $X$. Indeed, it follows from a result by Avramov (\cite{avr}) that if $\mathcal{X}$ is an algebraic $n$-stack locally of finite presentation over $k$\footnote{The result holds true if $\mathrm{Spec} \, k$ is replaced by any noetherian regular affine scheme, or more generally, if $\mathrm{Spec} \, k$ is replaced by any noetherian affine scheme $S$ and the morphism $\mathcal{X}\rightarrow S$ is assumed to be locally of finite flat dimension.}, and its stacky cotangent complex $\mathbb{L}$ is of \emph{bounded} amplitude, then $\mathbb{L}$ has in fact amplitude in $[-1,\infty)$ (actually in $[-1,n]$). The positive degrees are the (higher) stacky ones, while the negative degree $-1$ is there because it is already there e.g. for an lci scheme; the other negative degrees are forbidden by the boundedness assumption on $\mathbb{L}$ and Avramov's theorem.\\ Actually, \emph{any} derived extension $\mathcal{X}^{\mathrm{der}}$ of a stack $\mathcal{X}$ endows $\mathcal{X}$ with a deformation/obstruction theory (via the pullback of $\mathbb{L}_{\mathcal{X}^{\mathrm{der}}}$ along the closed immersion $\mathcal{X}\hookrightarrow\mathcal{X}^{\mathrm{der}}$); there is usually one natural such extension but there might be more, and different extensions could give a priori different deformation/obstruction theories. Conversely, from a $1$-stack endowed with a choice of a (perfect) obstruction theory one is able to essentially reconstruct the derived (quasi-smooth) extension whose derived cotangent complex induces the chosen obstruction theory. This might be compared, e.g., to the example of different obstruction theories (the \emph{standard} and the \emph{reduced} one) for the stack of stable maps to a $K3$-surface (\cite[2.2]{MP} and \cite[\S 4 and App. A]{mpt}), which should indeed correspond to different derived extensions of this stack. It would be interesting to unravel the moduli interpretation of the derived stack corresponding to the reduced obstruction theory. \\
We address the reader to \cite[4.1]{seattle} where most of the topics of this $\S$ are discussed more extensively.\\

\end{document}